%% file: baer-dahl.tex
\documentclass[10pt]{amsart}
\usepackage{a4,amssymb,amsfonts,pstricks}

\newcommand{\RR}{{\mathbb R}}

\renewcommand{\phi}{\varphi}
\newcommand{\eps}{\varepsilon}

\newcommand{\grad}{\operatorname{grad}}

\newcommand{\scal}{\operatorname{Scal}}

\newcommand{\<}{\left\langle}       
\renewcommand{\>}{\right\rangle}

\newcommand{\vol}{{\operatorname{vol}}}

\newtheorem{thm}{Theorem}
\newtheorem{lemma}[thm]{Lemma}

\newtheorem{remark}[thm]{Remark}
\newtheorem{remarks}[thm]{Remarks}
\newtheorem{definition}[thm]{Definition}
\newtheorem{notation}[thm]{Notation}
\newtheorem{example}[thm]{Example}

\begin{document}
%%%%%%%%%%%%%%%%%%%%%%%%%%%%%%%%%%%%%%%%%%%%%%%%%%%%%%%%%%%%%%%%%%%%%%%%%

\title
[The First Dirac Eigenvalue]
{The First Dirac Eigenvalue on Manifolds with Positive Scalar Curvature}

\author{Christian B{\"a}r}
\author{Mattias Dahl}

\address{
Universit{\"a}t Hamburg\\
FB Mathematik\\
Bundesstr.~55\\
20146 Hamburg\\
Germany
}
\address{
Institutionen f\"or Matematik\\
Kungl Tekniska H\"ogskolan\\
10044 Stockholm\\
Sweden
}
 
\email{baer@math.uni-hamburg.de} 
\email{dahl@math.kth.se}

\subjclass[2000]{53C27}

\keywords{Dirac operator, eigenvalue, positive scalar curvature, 
Friedrich's estimate}

\thanks{The first author has been partially supported by the Research and 
Training Networks HPRN-CT-2000-00101 ``EDGE'' and HPRN-CT-1999-00118 
``Geometric Analysis'' funded by the European Commission.}

\date{\today}

\begin{abstract}
We show that on every compact spin manifold admitting a Riemannian metric
of positive scalar curvature Friedrich's eigenvalue estimate for the Dirac 
operator can be made sharp up to an arbitrarily small given error by choosing 
the metric suitably.
\end{abstract}

\maketitle

\section{Introduction and Statement of the Result}

On an $n$-dimensional compact Riemannian spin manifold $M$ the Dirac operator 
has discrete real spectrum consisting only of eigenvalues of finite 
multiplicity.
If the manifold has positive scalar curvature, then $0$ lies in a spectral
gap, more precisely, Friedrich \cite{friedrich80a} showed that all 
eigenvalues $\lambda$ of the Dirac operator satisfy
$$
\lambda^2 \geq \frac{n}{4(n-1)}\min_M\scal .
$$
This inequality is sharp in the sense that there are examples in all dimensions
where equality is attained for the eigenvalue of smallest modulus.
The standard sphere provides such examples.
Equality in this estimate implies strong restrictions on the manifold.
The manifold must be Einstein and the corresponding eigenspinor then must 
be a Killing spinor.
The geometric types of manifolds admitting such Killing spinors are classified,
see \cite{baer93a}.
In particular, in even dimension $n\not= 6$ only the standard sphere carries
Killing spinors. 

The dimension dependent coefficient $c_n=\frac{n}{4(n-1)}$ in the estimate can
be improved if one imposes {\em geometric assumptions} on the metric.
Kirchberg \cite{kirchberg86a,kirchberg90a} showed that for K\"ahler
metrics $c_n$ can replaced by $\frac{n+2}{4n}$ if the complex dimension 
$\frac{n}{2}$ is odd, and by $\frac{n}{4(n-2)}$ if $\frac{n}{2}$ is even.
For quaternionic K\"ahler manifolds Kramer, Weingart, and Semmelmann
\cite{kramer-weingart-semmelmann99a} showed that $c_n$ can be improved to
$\frac{n+12}{4(n+8)}$.

Alexandrov, Grantcharov, and Ivanov \cite{alexandrov-grantcharov-ivanov98a}
showed that if there exists a parallel one-form on $M$, then $c_n$ can be
replaced by $c_{n-1}=\frac{n-2}{4(n-1)}$.
In a recent paper Moroianu and Ornea \cite{moroianu-ornea03a} weakened 
the assumption on the 1-form from parallel to harmonic with constant length.
They believed that the condition of having constant length is also
unnecessary and conjectured \cite[Conj.~1]{moroianu-ornea03a}:
All Dirac eigenvalues $\lambda$ on an $n$-dimensional compact Riemannian 
spin manifold {\em with nonzero first Betti-number} satisfy
$$
\lambda^2 \geq \frac{n-1}{4(n-2)}\min_M\scal .
$$
We show here that this conjecture is false.
More precisely, we prove

{\bf Theorem.}
{\em Let $M$ be a compact $n$-dimensional Riemannian spin manifold with 
positive scalar curvature.

Then there exists a smooth one-parameter family of Riemannian metrics
$g_\eps$ on $M$, $\eps\in(0,\eps_0]$, such that
\begin{itemize}
\item
$\scal_{g_\eps} \geq n(n-1)$
\item
$\frac{n^2}{4} \leq \lambda_1(D^2_{M,g_\eps}) \leq   \frac{n^2}{4} + \eps$
\end{itemize}
where $\lambda_1(D^2_{M,g_\eps})$ is the smallest eigenvalue of the 
square of the Dirac operator $D_{M,g_\eps}$ on $(M,g_\eps)$.
}

The lower eigenvalue bound is nothing but Friedrich's estimate.
Thus the theorem says that Friedrich's estimate can always be made
``almost sharp'' by choosing suitable metrics.
In particular, the dimension dependent coefficient $\frac{n}{4(n-1)}$ in
the estimate cannot be improved by imposing additional {\em topological 
assumptions} like bounds on the Betti numbers.

This is indeed remarkable since equality in Friedrich's estimate is so
restrictive on the manifold.
The manifold then has to be Einstein with positive Einstein
constant.
So it has positive Ricci curvature, hence its fundamental group is finite
by the Bonnet-Myers theorem.
In particular, the first Betti number must vanish.
On the other hand, our theorem says that ``almost equality'' in Friedrich's
estimate does not impose any topological restriction on the manifold.

\section{The Proof}

We start by proving two preliminary lemmas.

\begin{lemma}\label{lemmasubset}
Let $M$ be a compact Riemannian spin manifold of dimension $n\geq 3$,
let $p\in M$.

Then for each $\delta>0$ there exists $R_1(\delta)>0$ such that for each compact
$n$-dimensional Riemannian spin manifold $\widetilde{M}$ containing an 
isometric copy of $M-\bar{B}_r(p)$ as an open subset (with the same spin 
structure) where $B_r(p)$ is the geodesic ball about $p$ of radius $r$, 
$0<r\leq R_1(\delta)$, we have
$$
\lambda_1(D^2_{\widetilde{M}}) \leq \lambda_1(D^2_M) + \delta.
$$
Moreover, $R_1$ depends smoothly on $\delta$.
\end{lemma}

\begin{proof}
Let $\phi$ be a nontrivial eigenspinor on $M$ to the eigenvalue $\lambda$
of $D$ where $\lambda = \sqrt{\lambda_1(D^2_M)}$ or $\lambda = 
-\sqrt{\lambda_1(D^2_M)}$.
Fix $r_0>0$ such that $2r_0 < \mathrm{injrad}(p)$, the injectivity radius of 
$p$.
For $0<r \leq r_0$ choose a smooth cutoff-function $\chi \in C^\infty(M)$ 
satisfying
\begin{itemize}
\item
$0 \leq \chi \leq 1$ on all of $M$
\item
$\chi \equiv 0$ on $B_{r}(p)$
\item
$\chi \equiv 1$ on $M-B_{2r}(p)$
\item
$|\grad\chi| \leq 2/r$ on $M$
\end{itemize}
The spinor field $\chi\phi$ vanishes on $B_{r}(p)$ and can therefore also
be considered as a spinor on $\widetilde{M}$.
Therefore
\begin{eqnarray*}
\lambda_1(D^2_{\widetilde{M}}) &\leq&
\frac{\int_{\widetilde{M}}|D(\chi\phi)|^2\, dV}
{\int_{\widetilde{M}}|\chi\phi|^2\, dV} \\
&=&
\frac{\int_{{M}}|D(\chi\phi)|^2\, dV}
{\int_{{M}}|\chi\phi|^2\, dV} \\
&=&
\frac{\int_{{M}}|\chi D\phi + \grad\chi\cdot\phi|^2\, dV}
{\int_{{M}}\chi^2|\phi|^2\, dV} \\
&\stackrel{(*)}{=}&
\frac{\int_{{M}}\left(\chi^2\lambda^2|\phi|^2 + |\grad\chi|^2|\phi|^2\right)\, dV}
{\int_{{M}}\chi^2|\phi|^2\, dV} \\
&=&
\lambda_1(D^2_M) + 
\frac{\int_{{M}}|\grad\chi|^2|\phi|^2\, dV}{\int_{{M}}\chi^2|\phi|^2\,dV}\\
&\leq&
\lambda_1(D^2_M) + 
\frac{\frac{4}{r^2}\int_{B_{2r}(p)}|\phi|^2\, dV}
{\int_{M-B_{2r}(p)}|\phi|^2\,dV}\\
&\leq&
\lambda_1(D^2_M) + 
\frac{\frac{4}{r^2}\,\vol(B_{2r}(p))\,\|\phi\|^2_{L^\infty}}
{\int_{M-B_{2r_0}(p)}|\phi|^2\,dV}\\
&\leq&
\lambda_1(D^2_M) + C\cdot r^{n-2}
\end{eqnarray*}
where $C$ is a positive constant depending only on $M$, $r_0$, $p$, and $\phi$.
In $(*)$ we used that the mixed terms $\<\chi \lambda \phi, 
\grad\chi\cdot\phi\>$ and $\<\grad\chi\cdot\phi,\chi \lambda \phi\>$ 
cancel because Clifford multiplication with $\grad\chi$ is skew-symmetric.
The lemma follows.
\end{proof}

\begin{lemma}\label{lemmasphere}
Let $p\in S^n$, $n\geq 2$. 
Let $S>0$.

There exists a smooth one-parameter family of Riemannian metrics
$h_{\delta}$ on $S^n$, $\delta \in (0,\delta_0]$, such that
\begin{itemize}
\item  
$\scal_{h_{\delta}} \geq n(n-1)$
\item
$\scal_{h_{\delta}}|_{B_{R_2(\delta)}(p)} \geq S$ where $R_2(\delta)>0$
depends smoothly on $\delta$.
\item
$\lambda_1(D^2_{h_{\delta}}) \leq \frac{n^2}{4} + \delta$
\end{itemize}
\end{lemma}

\begin{proof}
The idea of the proof is this.
Consider the sphere with its standard metric embedded in Euclidean space
$\RR^{n+1}$, $S^n = \{(t,x_1,\ldots,x_n)\in\RR^{n+1}\,|\, |(t,x)|^2=1\}$.
Call $S^n_+ := \{(t,x_1,\ldots,x_n)\in\RR^{n+1}\,|\, |(t,x)|^2=1,\, t\geq 0\}$ 
the {\em northern hemisphere} and $S^n_- := 
\{(t,x_1,\ldots,x_n)\in\RR^{n+1}\,|\, |(t,x)|^2=1,\, t\leq 0\}$ the 
{\em southern hemisphere}.
Move the northern hemisphere by a small amount $\eta>0$ to the south 
(i.~e.\ in direction $-e_0 = (-1,0,\ldots,0)$) and cut off everything that 
gets moved south of the hyperplane $\{0\}\times\RR^n$.
This yields the northern cup $S^n_{+,\eta} := 
\{(t,x)\in\RR^{n+1}\,|\, |(t+\eta,x)|^2=1,\, t\geq -\eta\}$.
Similarly, move the southern hemisphere to the north by the amount $\eta$,
cut off everything north of the equatorial hyperplane, and obtain
$S^n_{-,\eta} := \{(t,x)\in\RR^{n+1}\,|\, |(t-\eta,x)|^2=1,\, t\leq \eta\}$.
The union $S^n_{+,\eta}\cup S^n_{-,\eta}$ is a hypersurface in 
$\RR^{n+1}$, singular along the equator $\partial S^n_{+,\eta} = \partial 
S^n_{-,\eta} = S^n_{+,\eta} \cap S^n_{-,\eta} $, smooth elsewhere with constant
Gauss curvature $K \equiv 1$ and constant mean curvature $H \equiv 1$.

\input{fig1}

Careful smoothing of this hypersurface in a neighborhood of the equator
yields a hypersurface $S^n_\eta$ diffeomorphic to $S^n$ such that
\begin{itemize}
\item
$\scal\geq n(n-1)$ everywhere
\item
$\scal\geq S$ in the $\eta^2(1-4\eta^2)$-tubular neighborhood of the equator
\item
$\frac{1}{\vol(S^n_\eta)}\int_{S^n_\eta} H^2\, dV \leq 1 + \mathrm{O}(\eta)$
\end{itemize}
If $p$ was placed on the equator in the beginning, which we may assume,
then the statements on the scalar curvature follow immediately.
The bound on the Dirac eigenvalue is a consequence of 
$$
\lambda_1(D^2_M) \leq \frac{n^2}{4\vol(M)}\int_M H^2\, dV
$$
for any compact oriented hypersurface $M$ of $\RR^{n+1}$, 
see \cite[Cor.~4.2]{baer98b}.

To make the smoothing of $S^n_{+,\eta}\cup S^n_{-,\eta}$ explicit we 
start with the case of dimension $n=2$.
The 2-sphere is a surface of revolution, parametrized by
$$
F : [-1,1]\times S^1 \to \RR^3, \quad
F(t,\theta) = (t,r(t)\cos(\theta),r(t)\sin(\theta))
$$
with $r(t) = \sqrt{1-t^2}$.
The general formula for the principal curvatures of a surface of revolution is 
$$ 
\kappa_t = \frac{-\ddot{r}}{(1+\dot{r}^2)^{3/2}}
$$
and 
$$
\kappa_\theta = \frac{1}{r\sqrt{1+\dot{r}^2}} .
$$
In our case of the sphere $\kappa_t = \kappa_\theta = 1$.
The smoothed hypersurface $S^2_\eta$ will also be a surface of revolution
with radius function $r_\eta:[-1+\eta,1-\eta] \to \RR$ chosen subject
to the following conditions:
\begin{enumerate}
\item
$r_\eta$ is smooth on $(-1+\eta,1-\eta)$ and even
\item
$r_\eta(t) = \sqrt{1-(t-\eta)^2}$ for $t\in [-1+\eta,-\eta]$
\item
$r_\eta(t) = \sqrt{1-(t+\eta)^2}$ for $t\in [\eta,1-\eta]$
\item
$\sqrt{1-4\eta^2} \leq r_\eta(t) \leq \frac{1}{\sqrt{1-4\eta^2}}$ 
for $t\in [-\eta,\eta]$
\item
$|\dot{r}_\eta(t)| \leq \frac{2\eta}{\sqrt{1-4\eta^2}}$ for $t\in [-\eta,\eta]$
\item
$\ddot{r}_\eta < 0$ on $(-1+\eta,1-\eta)$
\item
$\ddot{r}_\eta(t) = -2S$ for $t\in [-\eta^2,\eta^2]$
\item
$-2S \leq \ddot{r}_\eta(t) \leq \frac{-1}{(1-4\eta^2)^{3/2}}$ for 
$t\in [-\eta,\eta]$
\end{enumerate}

\input{fig2}

For (2) and (3) to make sense we assume that $\eta < \frac{1}{2}$.
By (2) we have $\dot{r}_\eta(-\eta) = \frac{2\eta}{\sqrt{1-4\eta^2}}$ and 
by (3) we get $\dot{r}_\eta(\eta) = \frac{-2\eta}{\sqrt{1-4\eta^2}}$.
Thus (6) implies (5).

Moreover, by (2) and (3) we obtain $r_\eta(\pm\eta) = \sqrt{1-4\eta^2}$
and so, by (6), we see $r_\eta \geq \sqrt{1-4\eta^2}$ on $[-\eta,\eta]$.
By (5) we see that $r_\eta(t) \leq r_\eta(-\eta) + \eta\frac{2\eta}
{\sqrt{1-4\eta^2}} \leq \frac{1}{\sqrt{1-4\eta^2}}$ for $t\in [-\eta,0]$
and similarly for $t\in [0,\eta]$.
Hence (4) is also a consequence of (2), (3), and (6).

Condition (6) follows from (2), (3), and (8).

To make sure that conditions (2), (3), (7), and (8) can be realized we need
to assume $S > 1$, which we do without loss of generality, and that
$\eta$ is so small that
\begin{itemize}
\item[(a)]
$2S\eta^2 < -\dot{r}_\eta(\eta) = \frac{2\eta}{\sqrt{1-4\eta^2}}$
\item[(b)]
$-2S\eta < \dot{r}_\eta(\eta)$, i.~e.\ $S > \frac{1}{\sqrt{1-4\eta^2}}$
\item[(c)]
$-2S\eta^2 > \dot{r}_\eta(\eta) - \ddot{r}_\eta(\eta)\cdot (\eta-\eta^2)$,
i.~e. $S<\frac{\frac1\eta + 1 - 8\eta}{2(1-4\eta^2)^{3/2}}$
\end{itemize}

\input{fig3}

Conditions (1) - (8) together with the explicit formulas for the 
principal curvatures imply
$$
\begin{array}{cl}
1 \leq \kappa_t(t,\theta) \leq 2S & \mbox{for }t \in [-\eta,\eta]\\
\kappa_t(t,\theta) \geq 2S(1-4\eta^2)^{3/2} & 
\mbox{for }t \in [-\eta^2,\eta^2]\\
1-4\eta^2 \leq \kappa_\theta(t,\theta) \leq \frac{1}{\sqrt{1-4\eta^2}}&
\mbox{for }t \in [-\eta,\eta]
\end{array}
$$
At this point we return to the general case of dimension $n\geq 2$.
If $n\geq 3$, then we perform exactly the same smoothing, i.~e.\ we
use the same warping function $r_\eta(t)$ where $t$ is the first Cartesian
coordinate in $\RR^{n+1}$.
For every 3-dimensional vector subspace $E\subset\RR^{n+1}$ containing $e_0$ 
the reflection across $E$ is an isometry of $\RR^{n+1}$ leaving the smoothed 
hypersurface $S^n_\eta$ invariant.
Hence the fixed point set $E\cap S^n_\eta$ is a totally geodesic submanifold
of $S^n_\eta$ so that its principal curvatures are also principal curvatures 
of $S^n_\eta$.
Therefore the principal curvatures $\kappa_1,\ldots,\kappa_n$ of $S^n_\eta$ 
satisfy
$$
\begin{array}{cl}
1 \leq \kappa_1(t,\theta) \leq 2S & \mbox{for }t \in [-\eta,\eta]\\
\kappa_1(t,\theta) \geq 2S(1-4\eta^2)^{3/2} & 
\mbox{for }t \in [-\eta^2,\eta^2]\\
1-4\eta^2 \leq \kappa_j(t,\theta) \leq \frac{1}{\sqrt{1-4\eta^2}}&
\mbox{for }t \in [-\eta,\eta]
\end{array}
$$
$j=2,\ldots,n$.
This implies for the scalar curvature $\scal = \sum_{i\not= j}\kappa_i\kappa_j$
and the mean curvature $H = \frac1n\sum_j\kappa_j$
$$
\begin{array}{cl}
\scal \geq n(n-1)(1-4\eta^2)^2 & \mbox{for }t \in [-\eta,\eta]\\
\scal > S & \mbox{for }t \in [-\eta^2,\eta^2]\\
H \leq 2S & \mbox{for }t \in [-\eta,\eta]\\
H = 1 &  \mbox{for }t \in [-1+\eta,-\eta] \cup [\eta,1-\eta]\\
\scal = n(n-1) &  \mbox{for }t \in [-1+\eta,-\eta] \cup [\eta,1-\eta]
\end{array}
$$
The volume of the part of $S^n_\eta$ contained in the strip $-\eta \leq t
\leq \eta$ is of order $\mathrm{O}(\eta)$.
Thus 
\begin{eqnarray*}
\frac{1}{\vol{S^n_\eta}}\int_{S^n_\eta}H^2 dV 
&=&
1 + \frac{1}{\vol{S^n_\eta}}\int_{S^n_\eta}(H^2-1)\, dV \\
&\leq&
1 + \frac{1}{\vol{S^n_\eta}}\int_{S^n_\eta\cap\{-\eta \leq t\leq \eta\}}
2S\, dV \\
&=&
1 + \mathrm{O}(\eta) .
\end{eqnarray*}
Hence $S^n_\eta$ satisfies
\begin{itemize}
\item
$\scal\geq n(n-1)(1-4\eta^2)^2$ everywhere
\item
$\scal\geq S$ in the $\eta^2$-tubular neighborhood of the equator
\item
$\frac{1}{\vol(S^n_\eta)}\int_{S^n_\eta} H^2\, dV \leq 1 + \mathrm{O}(\eta)$
\end{itemize}
Multiplying $S^n_\eta$  by the factor $(1-4\eta^2)$ yields
\begin{itemize}
\item
$\scal\geq n(n-1)$ everywhere
\item
$\scal\geq \frac{S}{(1-4\eta^2)^2}\geq S$ in the $\eta^2(1-4\eta^2)$-tubular 
neighborhood of the equator
\item
$\frac{1}{\vol(S^n_\eta)}\int_{S^n_\eta} H^2\, dV \leq 
\frac{1 + \mathrm{O}(\eta)}{(1-4\eta^2)^2} = 1 + \mathrm{O}(\eta)$
\end{itemize}
Substituting a suitable multiple of $\eta$ by $\delta$ concludes the proof.
\end{proof}

{\em Proof of the Theorem.}
Without loss of generality we may assume that $M$ is connected.
If $M$ is 2-dimensional it must be the 2-sphere.
The constant family $g_\eps = g_{\mathrm{can}}$ does the job where 
$g_{\mathrm{can}}$ is the standard metric of constant Gauss curvature 1 because
\begin{itemize}
\item
$\scal_{g_{\mathrm{can}}} \equiv 2 = n(n-1)$
\item
$\lambda_1(D^2_{S^2,g_{\mathrm{can}}}) = 1 = \frac{n^2}{4}$
\end{itemize}

From now on let $n\geq 3$.
Pick a Riemannian metric $h$ on $M$ such that $\scal_h \geq 2n(n-1)$.
Apply Lemma~\ref{lemmasphere} with $S=2n(n-1)$ and obtain a smooth
one-parameter family of metrics $h_\delta$ on $S^n$ such that
\begin{itemize}
\item  
$\scal_{h_{\delta}} \geq n(n-1)$
\item
$\scal_{h_{\delta}}|_{B_{R_2(\delta)}(p)} \geq 2n(n-1)$
\item
$\lambda_1(D^2_{h_{\delta}}) \leq \frac{n^2}{4} + \delta$
\end{itemize}
Let $r(\delta)$ be a smooth function of $\delta$ such that $0<r(\delta) \leq
\min\{R_1(\delta),R_2(\delta)\}$.
Now form the connected sum $\widetilde{M}$ of $(S^n,h_\delta)$ and $(M,h)$ such
that the metric $h_\delta$ remains unchanged outside $B_{r(\delta)}(p)$.
Thus $\widetilde{M}$ contains an isometric copy of $S^n-\bar{B}_{r(\delta)}(p)$
and by Lemma~\ref{lemmasubset} we obtain
$$
\lambda_1(D^2_{\widetilde{M}})  \leq \lambda_1(D^2_{S^n,h_\delta}) + \delta
\leq \frac{n^2}{4} + 2\delta .
$$
Performing connected sums in dimension $n\geq 3$ (or, more generally,
surgery in codimension $\leq 3$) does not decrease scalar curvature
too much if the metric on the connected sum is chosen carefully, see 
\cite[Proof of Theorem~A]{gromov-lawson80a} and
\cite[Proof of Theorem~3.1]{rosenberg-stolz98a}.
Since $\scal \geq 2n(n-1)$ on $B_{r(\delta)}(p)$ as well as on $M$
we may assume that the metric on $\widetilde{M}$ still has scalar
curvature $\scal \geq n(n-1)$.
Moreover, the construction of the metric on $\widetilde{M}$ can be done
smoothly in $\delta$.
Since $\widetilde{M}$ is diffeomorphic to the original manifold $M$
the substitution $\eps=2\delta$ yields the theorem.
$\hfill\Box$

%%%%%%%%%%%%%%%%%%%%%%%%%%%%%%%%%%%%%%%%%%%%%%%%%%%%%%%%%%%%%%%%%%%%%%%%%
%\bibliographystyle{amsplain}
%\bibliography{meine.bib,allg.bib}  

\providecommand{\bysame}{\leavevmode\hbox to3em{\hrulefill}\thinspace}

%%%%%%%%%%%%%%%%%%%%%%%%%%%%%%%%%%%%%%%%%%%%%%%%%%%%%%%%%%%%%%%%%%%%%%%%%

\end{document}

%% file: fig1.tex
\begin{center}
\begin{pspicture}(-5,-1.8)(5,2)

\psset{unit=1cm}
%\psgrid(-5,-4)(5,4)

\psellipse(0,0)(2,0.5)
\psframe[linecolor=white,fillstyle=solid,fillcolor=white](-2.1,1.1)(2.1,0)
\psellipse[linestyle=dotted](0,0)(2,0.5)
\psarc(0,-0.54){2.08}{15}{165}
\psarc(0,0.54){2.08}{195}{-15}

\rput(2,1){$S^n_{+,\eta}$}
\rput(2,-1.1){$S^n_{-,\eta}$}

\end{pspicture}
\end{center}

\begin{center}
{\sc Fig.~1}
\end{center}

%% file: fig2.tex
\begin{center}
\begin{pspicture}(-5,-0.5)(5,4.5)

\psset{unit=2cm}
%\psgrid(-5,-1.8)(5,2)

\psline{->}(-2.5,0)(2.5,0)
\psline(0,0)(0,2.2)
\psarc[linecolor=lightgray](0.5,0){2}{94}{180}
\psarc(0.5,0){2}{115}{180}
\psarc[linecolor=lightgray](-0.5,0){2}{0}{86}
\psline[linestyle=dotted](-0.35,0)(-0.35,2)
\psarc(-0.5,0){2}{0}{65}
\psline[linestyle=dotted](0.35,0)(0.35,2)
\psecurve(-1,1)(-0.35,1.81)(0.35,1.81)(1,1)
\psline[linestyle=dotted](1.5,1.8)(0.35,1.8)
\psline[linestyle=dotted](0,2)(1.4,2)

\rput(0,-0.2){$0$}
\rput(-0.35,-0.2){$-\eta$}
\rput(0.35,-0.2){$\eta$}
\rput(-1.6,-0.2){$-1+\eta$}
\rput(1.5,-0.2){$1-\eta$}
\rput(-1.4,1){$r_\eta$}
\rput(1.9,1.8){$\sqrt{1-4\eta^2}$}
\rput(1.5,2){$1$}
\rput(2.5,0.2){$t$}

\end{pspicture}
\end{center}

\begin{center}
Graph of ${r}_\eta$
\end{center}
\begin{center}
{\sc Fig.~2}
\end{center}

%% file: fig3.tex
\begin{center}
\begin{pspicture}(-5,-3)(5,3)

\psset{unit=2cm}
%\psgrid(-5,-0.5)(5,4.5)

\psline{->}(-2.5,0)(2.5,0)
\psline(0,-1.5)(0,1.5)
\psline[linecolor=lightgray](-1.5,1.5)(1.5,-1.5)
\psline(-0.5,0.5)(0.5,-0.5)
\psecurve(0,-1)(1.5,-1)(2,-1.2)(3,-3)
\psline[linecolor=lightgray](1.5,-1)(0.5,-0.8)
\psecurve(0,0.2)(0.5,-0.5)(1.5,-1)(2,-1)
\psdot(0.5,-0.5)
\psdot(0.5,-0.8)
\psdot(1.5,-1)
\psdot(1.5,-1.5)
\psecurve(0,1)(-1.5,1)(-2,1.2)(-3,3)
\psline[linecolor=lightgray](-1.5,1)(-0.5,0.8)
\psecurve(0,-0.2)(-0.5,0.5)(-1.5,1)(-2,1)

\psline[linestyle=dotted](0.5,-0.8)(0.5,0)
\psline[linestyle=dotted](1.5,-1)(1.5,0)
\rput(0.5,0.2){$\eta^2$}
\rput(1.5,0.2){$\eta$}
\rput(2.5,0.2){$t$}

\psline[linestyle=dotted](0.5,-0.5)(-0.5,-0.5)
\rput(-0.8,-0.5){$-2S\eta^2$}
\psline[linestyle=dotted](0.5,-0.8)(-1,-0.8)
\rput(-1.95,-0.8){$\dot{r}_\eta(\eta) - \ddot{r}_\eta(\eta)\cdot (\eta-\eta^2)$}
\psline[linestyle=dotted](1.5,-1)(2.1,-1)
\rput(2.35,-1){$\dot{r}_\eta(\eta)$}
\psline[linestyle=dotted](1.5,-1.5)(2.1,-1.5)
\rput(2.35,-1.5){$-2S\eta$}

\rput(-1.7,0.9){$\dot{r}_\eta$}

\end{pspicture}
\end{center}

\begin{center}
Graph of $\dot{r}_\eta$
\end{center}

\begin{center}
{\sc Fig.~3}
\end{center}